# Elementary properties of power series fields over finite fields

Franz-Viktor Kuhlmann

2. 6. 1998

## Abstract

In spite of the analogies between $\mathbb{Q}_p$ and $\mathbb{F}_p((t))$ which became evident through the work of Ax and Kochen, an adaptation of the complete recursive axiom system given by them for $\mathbb{Q}_p$ to the case of $\mathbb{F}_p((t))$ does not render a complete axiom system. We show the independence of elementary properties which express the action of additive polynomials as maps on $\mathbb{F}_p((t))$. We formulate an elementary property expressing this action and show that it holds for all maximal valued fields. We also discuss the action of arbitrary polynomials on valued fields.

## Contents

## 1  Elementary properties and additive polynomials

In this paper, we work with valued fields $(K, v)$, denoting the value group by $vK$, the residue field by $Kv$ and the valuation ring by $\mathcal{O}_v$ or just $\mathcal{O}$. For elements $a \in K$, the value is denoted by $va$, and the residue by $av$. We will use the classical additive (Krull) way of writing valuations. That is, the value group is an additively written ordered abelian group, the homomorphism property of $v$ reads as $vab = va + vb$, and the ultrametric triangle law reads as $v(a+b) \geq \min\{va, vb\}$. Further, we have the rule $va = \infty \Leftrightarrow a = 0$. We fix a language $\mathcal{L}$ of valued fields (or valued rings) which contains a relation symbol $\mathcal{O}(X/Y)$ for valuation divisibility. That is, $\mathcal{O}(a/b)$ will say that $va \geq vb$, or equivalently, that $a/b$ is an element of the valuation ring. We will write $\mathcal{O}(X)$ in the place of $\mathcal{O}(X/1)$ (note that $\mathcal{O}(a/1)$ says that $va \geq v1 = 0$, i.e., $a \in \mathcal{O}_v$).

Let $\mathbb{F}_p$ denote the field with $p$ elements. The power series field $\mathbb{F}_p((t))$, also called "field of formal Laurent series over $\mathbb{F}_p$", carries a canonical valuation $v_t$, the $t$-adic valuation (we write $v_t t = 1$). $(\mathbb{F}_p((t)), v_t)$ is a complete discretely valued field, with value group $v_t \mathbb{F}_p((t)) = \mathbb{Z}$ (that is what "discretely valued" means) and residue field $\mathbb{F}_p((t))v_t = \mathbb{F}_p$. At the first glimpse, such fields may appear to be the best known objects in valuation theory. Nevertheless, the following prominent questions about the elementary theory $\mathrm{Th}(\mathbb{F}_p((t)), v_t)$ are still unanswered:

**Is $\mathrm{Th}(\mathbb{F}_p((t)), v_t)$ decidable? Is it model complete? Does $(\mathbb{F}_p((t)), v_t)$ admit quantifier elimination in $\mathcal{L}$ or in a natural extension of $\mathcal{L}$? Does there exist an elementary class of valued fields, containing $(\mathbb{F}_p((t)), v_t)$ and satisfying some Ax–Kochen–Ershov principle?**

By an Ax–Kochen–Ershov principle for a class $\mathcal{K}$ of valued fields we mean a principle of the form

$$(K, v), (L, v) \in \mathcal{K} \text{ with } vK \equiv vL,\ Kv \equiv Lv \text{ implies that } (K, v) \equiv (L, v)$$

or a similar version with $\prec$ or $\prec_\exists$ ("existentially closed in") in the place of $\equiv$. Here, $vK$ denotes the value group of $(K, v)$, and the language is that of ordered groups. Further, $Kv$ denotes the residue field of $(K, v)$, and the language is that of rings or of fields. For example, the elementary class of henselian fields with residue fields of characteristic 0 satisfies all of these Ax–Kochen–Ershov principles (cf. [AK], [E], [KP], [K2]).

Encouraged by the similarities between $\mathbb{F}_p((t))$ and the field $\mathbb{Q}_p$ of $p$-adics, one might try to give a complete axiomatization for $\mathrm{Th}(\mathbb{F}_p((t)), v_t)$ by adapting the well known axioms for $\mathrm{Th}(\mathbb{Q}_p, v_p)$. They express that $(\mathbb{Q}_p, v_p)$ is a henselian valued field of characteristic 0 with value group a $\mathbb{Z}$-group (i.e., an ordered abelian group elementarily equivalent to $\mathbb{Z}$), and residue field $\mathbb{F}_p$. They also express that $vp = 1$ (the smallest positive element in the value group). This is not relevant for $\mathbb{F}_p((t))$ since there, $p \cdot 1 = 0$. Nevertheless, we may add a constant name $t$ to $\mathcal{L}$ so that one can express by an elementary sentence that $vt = 1$.

A naive adaptation would just replace "characteristic 0" by "characteristic $p$" and $p$ by $t$. But there is an elementary property of valued fields that is satisfied by all valued fields of residue characteristic 0 and all formally $p$-adic fields, but not by all valued fields in general. It is the property of being **defectless**. A valued field $(K, v)$ is called defectless if the **fundamental equality**

$$n = \sum_{i=1}^{\mathrm{g}} \mathrm{e}_i \mathrm{f}_i$$

holds for every finite extension $L|K$, where $n = [L : K]$ is the degree of the extension, $v_1, \ldots, v_\mathrm{g}$ are the distinct extensions of $v$ from $K$ to $L$, $\mathrm{e}_i = (v_i L : vK)$ are the respective ramification indices, and $\mathrm{f}_i = [Lv_i : Kv]$ are the respective inertia degrees. (Note that $\mathrm{g} = 1$ if $(K, v)$ is henselian.) There is a simple example, probably already due to F. K. Schmidt, which shows that there are henselian discretely valued fields of positive characteristic which are not defectless.

However, each power series field with its canonical valuation is henselian and defectless. In particular, $(\mathbb{F}_p((t)), v_t)$ is defectless. For a less naive adaptation of the axiom system of $\mathbb{Q}_p$, we will thus add "defectless". We obtain the following axiom system in the language $\mathcal{L}(t)$:

$$\left.\begin{array}{l} (K,v) \text{ is a henselian defectless valued field} \\ K \text{ is of characteristic } p \\ vK \text{ is a } \mathbb{Z}\text{-group} \\ Kv = \mathbb{F}_p \\ vt \text{ is the smallest positive element in } vK \ . \end{array}\right\} \tag{1}$$

Let us note that also $(\mathbb{F}_p(t), v_t)^h$, the henselization of $(\mathbb{F}_p(t), v_t)$, satisfies these axioms. It was common knowledge since some time that this is a defectless field, and the proof of this fact is not all too hard. But it can also be deduced from a more general principal, the "generalized Grauert–Remmert Stability Theorem" (see [K2] for this theorem and its proof, and [K5], [K7] for further applications). It is also well-known that $(\mathbb{F}_p(t), v_t)^h$ is existentially closed in $(\mathbb{F}_p((t)), v_t)$; for an easy proof see [K2]. But it is not known whether $(\mathbb{F}_p((t)), v_t)$ is an elementary extension of $(\mathbb{F}_p(t), v_t)^h$.

In fact, it did not seem unlikely that axiom system (1) could be complete, until we proved in [K1]:

**Theorem 1** *The axiom system (1) is not complete.*

We wish to show how this result is obtained and which additional previously unknown elementary properties of $\mathbb{F}_p((t))$ have been discovered.

We start by noting that for $K = \mathbb{F}_p((t))$, the elements $1, t, t^2, \ldots, t^{p-1}$ form a basis of the field extension $K|K^p$. Thus,

$$K = K^p \oplus tK^p \oplus \ldots \oplus t^{p-1}K^p \ . \tag{2}$$

It follows that the $\mathcal{L}(t)$-sentence

$$\forall X \exists X_0 \ldots \exists X_{p-1} \ \ X = X_0^p + tX_1^p + \ldots + t^{p-1}X_{p-1}^p \tag{3}$$

holds in $K$.

Since the Frobenius $x \mapsto x^p$ is an endomorphism of every field $K$ of characteristic $p$, it follows that for every $i$ the polynomial $t^i X^p$ is **additive**. A polynomial $f(X) \in K[X]$ is called additive if $f(a+b) = f(a) + f(b)$ for all $a, b$ in any extension field of $K$. The additive polynomials in $K[X]$ are precisely the polynomials of the form

$$\sum_{i=0}^m c_i X^{p^i} \quad \text{with } c_i \in K, \ m \in \mathbb{N}$$

(cf. [L], VIII, §11). If $K$ is infinite, then $f(X) \in K[X]$ is additive if and only if $f(a+b) = f(a) + f(b)$ for all $a, b \in K$. For further details about additive polynomials, see [O], [W1], [W2] and [K2].

Now it is a natural question to ask what might happen if we replace the polynomials $t^i X^p$ in (3) by other additive polynomials. Apart from the additive polynomials $cX^{p^n}$, the most important is the **Artin-Schreier polynomial** $\wp(X) := X^p - X$. Lou van den Dries observed that if $k$ is a field of characteristic $p$ such that $\wp(k) := \{\wp(x) \mid x \in k\} = k$, then the $\mathcal{L}(t)$-sentence

$$\forall X \exists X_0 \ldots \exists X_{p-1} \ \ X = X_0^p - X_0 + tX_1^p + \ldots + t^{p-1}X_{p-1}^p \tag{4}$$

holds in $k((t))$. However, he found that he was not able to eliminate the quantifiers in this assertion (with respect to a version of axiom system (1) where "$Kv = \mathbb{F}_p$" is replaced by "$Kv \equiv k$"). Observe that $\wp(\mathbb{F}_p) = \{0\} \neq \mathbb{F}_p$. To get an assertion valid in $\mathbb{F}_p((t))$, we have to introduce a corrective summand $Y$:

$$\forall X \exists Y \exists X_0 \ldots \exists X_{p-1} \ X = Y + X_0^p - X_0 + tX_1^p + \ldots + t^{p-1}X_{p-1}^p \wedge \mathcal{O}(Y) \tag{5}$$

**Lemma 2** *The $\mathcal{L}(t)$-sentence (5) holds for every intermediate field $(K, v)$ between the fields $(\mathbb{F}_p(t), v_t)$ and $(\mathbb{F}_p((t)), v_t)$.*

Proof: Take $x \in K$. If $vx \geq 0$, then we set $y = x$ and $x_i = 0$ to obtain that $x = y = y + x_0^p - x_0 + tx_1^p + \ldots + t^{p-1}x_{p-1}^p$ with $vy \geq 0$. For $vx < 0$, we can proceed by induction on $-vx$ since $vK = \mathbb{Z}$. Suppose that $m \in \mathbb{N}$ and that we have shown the assertion to hold for every $x$ of value $vx > -m$. Take $x \in K$ such that $vx = -m$. There is $\ell \in \{0, \ldots, p-1\}$ such that $vx \equiv \ell$ modulo $p\mathbb{Z} = pvK$. Choose some $z \in K$ such that $vx = \ell + pvz = vt^\ell z^p$. Then $v(x/t^\ell z^p) = 0$, and the residue of $x/t^\ell z^p$ is some element $j \in \mathbb{F}_p$. It follows that $v(x/t^\ell(jz)^p - 1) = v(j^{-1}x/t^\ell z^p - 1) > 0$. Hence, $v(x - t^\ell(jz)^p) > vt^\ell(jz)^p = vx$. If $\ell > 0$, then we set $x' := x - t^\ell(jz)^p$, so that $vx' > vx$. If $\ell = 0$, then we set $x' := x - (jz)^p + jz$; since $vjz < 0$, we have that $vx = v(jz)^p < vjz$ and thus again, $vx' \geq \min\{v(x - (jz)^p), vjz\} > vx$. So by induction hypothesis, there are $y, x_0' \ldots x_{p-1}'$ such that $vy \geq 0$ and $x' = y + (x_0')^p - x_0' + t(x_1')^p + \ldots + t^{p-1}(x_{p-1}')^p$. We set $x_\ell = x_\ell' + z$ and $x_i = x_i'$ for $i \neq \ell$, to obtain by additivity that $x = y + x_0^p - x_0 + tx_1^p + \ldots + t^{p-1}x_{p-1}^p$. $\square$

This lemma shows that in analogy to (2), every intermediate field $(K, v)$ between $(\mathbb{F}_p(t), v_t)$ and $(\mathbb{F}_p((t)), v_t)$ satisfies:

$$K = \mathcal{O} + \wp(K) + tK^p + \ldots + t^{p-1}K^p. \tag{6}$$

If in addition $(K, v)$ is henselian, then we can improve this representation to

$$K = \mathbb{F}_p + \wp(K) + tK^p + \ldots + t^{p-1}K^p. \tag{7}$$

This is seen as follows. Using Hensel's Lemma, one proves that the valuation ideal $\mathcal{M}$ of any henselian field $(K, v)$ is contained in $\wp(K)$. On the other hand, $Kv = \mathbb{F}_p$ implies that $\mathcal{O} = \mathbb{F}_p + \mathcal{M}$. Consequently, $\mathbb{F}_p + \wp(K) = \mathcal{O} + \wp(K)$.

Theorem 1 is proved by constructing a valued field $(L, v)$ which satisfies axiom system (1) but not sentence (5):

**Theorem 3** *Take $(K, v)$ to be $(\mathbb{F}_p(t), v_t)^h$ or $(\mathbb{F}_p((t)), v_t)$. Then there exists an extension $(L, v)$ of $(K, v)$ such that:*
*a) $L|K$ is a regular extension of transcendence degree 1,*
*b) $1, t, t^2, \ldots, t^{p-1}$ is a basis of $L|L^p$,*
*c) $(L, v)$ is henselian defectless,*
*d) the value group $vL$ a $\mathbb{Z}$-group,*
*e) the residue field $Lv$ is again equal to $\mathbb{F}_p$,*
*f) sentence (5) does not hold in $(L, v)$.*

We have chosen to construct this extension also over $(\mathbb{F}_p(t), v_t)^h$ because this leads to a quite small valued field, having only transcendence degree 2 over its prime field. This allows us to apply it also to the problem of local uniformization in positive characteristic (cf. [K3] and [K5]). Note that a field extension $L|K$ is said to be **regular** if it is linearly disjoint from the algebraic closure of $K$, that is, if it is separable and $K$ is relatively algebraically closed in $L$.

We will give the construction of $(L, v)$ in Chapter 4; it is taken over from [K1]. The basic idea is to start with a simple transcendental extension $K(x)|K$ and extend the valuation $v$ such that $vx > vK$, so that $vK(x)$ is the lexicographically ordered product $\mathbb{Z} \times \mathbb{Z}$. Then automatically, $K(x)v = Kv = \mathbb{F}_p$. Passing to the henselization of $(K(x), v)$ doesn't change the value group and residue field. By adjoining $n$-th roots of suitable elements with value not in $vK$, we enlarge the value group without changing the residue field. While adjoining $p^m$-th roots, we build in a "twist" which in the end guarantees that $x$ cannot be of the form as stated in assertion (5). A major problem in the construction is how to obtain that the constructed field is defectless. (It is henselian, being an algebraic extension of the henselization of $(K(x), v)$.) To solve this problem, we use a characterization of defectless valued fields of positive characteristic, which we derived in [K1]. It is based on a classification of proper finite **immediate** extensions of henselian fields; an extension of valued fields is called immediate if it leaves value group and residue field unchanged. Such extensions violate the fundamental equality in the worst possible way, since $n > 1$ while e = f = g = 1.

Since $(L, v)$ does not satisfy (5), it cannot be an elementary extension of $(\mathbb{F}_p((t)), v_t)$. This contrasts the fact that, according to another theorem proved in [K1], $(\mathbb{F}_p((t)), v_t)$ is existentially closed in $(L, v)$.

Having seen that the sentence (5) is independent of the axioms in (1), we now pursue two main questions. The first of them is:

**A)** Are there further assertions similar to (5) and independent of (1)? What happens if we replace the additive polynomials $\wp(X)$, $tX^p$, ..., $t^{p-1}X^p$ appearing in (5) by other additive polynomials? Which corrective summands are then needed? Can we find a form that asserts essentially the same but dispenses with the use of the corrective summands $Y$, $\mathcal{O}$, $\mathbb{F}_p$ in (5), (6) and (7)?

Before we formulate the second question, let us give some background. In the model theory of valued fields, the **maximal fields** play a crucial role. These are valued fields not admitting any proper immediate extensions. It was shown by Krull [KR] that every valued field has at least one maximal immediate extension; this must be a maximal field. (Later, Gravett [G] gave a beautiful short proof replacing Krull's complicated argument.) As it is the case for power series fields (which in fact are maximal), also all maximal fields are henselian defectless.

A valued field $(K, v)$ is called **algebraically maximal** if it admits no proper immediate algebraic extension. As the henselization is an immediate algebraic extension, every algebraically maximal field is henselian. On the other hand, every henselian defectless field is algebraically maximal since every finite immediate extension would satisfy $e = f = g = 1$. But F. Delon [D] gave an example of an algebraically maximal field which is not defectless (this example can also be found in [K2]). Certainly, an algebraically maximal field is not necessarily maximal.

For the elementary classes of (see [K2] for missing definitions)

$$\left.\begin{array}{l}\bullet \text{ all henselian fields with residue characteristic 0 (cf. [AK], [E], [KP]),}\\ \bullet \text{ all henselian formally } p\text{-adic fields (cf. [AK], [E]),}\\ \bullet \text{ all henselian finitely ramified fields (cf. [E], [Z]),}\\ \bullet \text{ all algebraically maximal Kaplansky fields (cf. [E], [Z]),}\end{array}\right\} \quad (8)$$

most proofs of their good model theoretical properties work, implicitly or explicitly, with the following fact:

*the maximal immediate extensions of fields with residue characteristic 0, formally p-adic fields, finitely ramified fields and Kaplansky fields are unique up to isomorphism.*

(This actually follows from the fact that for such fields the maximal immediate algebraic extensions are unique up to isomorphism; cf. [KPR].) But uniqueness of maximal immediate extensions does not hold for arbitrary valued fields. In fact, there exist henselian fields with value group a $\mathbb{Z}$-group and residue field $\mathbb{F}_p$ which admit infinitely many non-isomorphic maximal immediate extensions.

Because of this special role of maximal fields, it would be important to know whether all maximal fields satisfy assertions similar to (5). But for an arbitrary maximal field $(M, v)$, also the $p$-degree $[M : M^p]$ is arbitrary and thus, the basis $1, t, \ldots, t^{p-1}$ has to be replaced adequately. On the other hand, elementary properties like "henselian" and "defectless" hold simultaneously for all maximal fields (see [K2] for the proofs), and they can be formulated without referring to the $p$-degree. So we ask:

**B)** Do all maximal fields satisfy assertions similar to (5)? Is there a way to formulate these assertions simultaneously for all maximal fields, not involving the $p$-degree?

In order to formulate our answer to these questions, we have to introduce some notation. Take a valued field $(K, v)$ of characteristic $p > 0$ and additive polynomials $f_0, \ldots, f_n \in K[X]$. We define an $\mathcal{L}$-formula

$$\mathrm{pd}(z_0, \ldots, z_n, z'_0, \ldots, z'_n) \;:\Leftrightarrow\; v\left(\sum_{i=0}^n z_i - \sum_{i=0}^n z'_i\right) > v\sum_{i=0}^n z_i \;\wedge\; v\sum_{i=0}^n z'_i = \min_i v z'_i$$

and an $\mathcal{L}(K)$-sentence

$$\mathrm{PD}(f_0, \ldots, f_n) \;:\Leftrightarrow\; \forall X_0, \ldots, X_n \exists Y_0, \ldots, Y_n \; \mathrm{pd}(f_0(X_0), \ldots, f_n(X_n), f_0(Y_0), \ldots, f_n(Y_n)) \;.$$

To understand the meaning of PD observe that $v \sum_{i=0}^n z_i \geq \min_i v z_i$ by the ultrametric triangle law, but that equality need not hold in general. In this situation, we would like to replace the $z_i$'s by $z'_i$'s such that $\sum_{i=0}^n z_i = \sum_{i=0}^n z'_i$ and $v \sum_{i=0}^n z'_i = \min_i v z'_i$. If one restricts the choice of the $z'_i$'s to certain sets (e.g., the images of the $f_i$'s), then this might not always be possible. Asking for the equality of the sums is quite strong; for our purposes, a weaker condition will suffice. We replace the equality by the expression $v(\sum_{i=0}^n z_i - \sum_{i=0}^n z'_i) > v \sum_{i=0}^n z_i$. This means that the new sum "approximates" the old, in a certain sense. Note that this implies that $v \sum_{i=0}^n z_i = v \sum_{i=0}^n z'_i$.

At this point, observe that the images $f_i(K)$ of $K$ under $f_i$ are subgroups of the additive group of $K$ because the $f_i$'s are additive. Now if we have subgroups $G_0, \ldots, G_n$ then we call their sum **direct (as valued groups)** if $v \sum_{i=0}^n z_i = \min_i v z_i$ for every choice of $z_i \in G_i$. In fact, $K = \mathbb{F}_p((t))$ is the direct sum of the subgroups $K^p, tK^p, \ldots, t^{p-1}K^p$ not only in the ordinary sense, but also as valued groups (see Lemma 16 in Section 3). On the other hand, the sum of the subgroups $\wp(K), tK^p, \ldots, t^{p-1}K^p$ is not direct since $t^i \mathcal{O} \subset \mathcal{M} \subset \wp(K)$ for all $i \geq 1$. Therefore, we introduce the notion **pseudo direct**: we call the sum of the $G_i$ pseudo direct if for every choice of $z_i \in G_i$ there are $z_i' \in G_i$ such that $\mathrm{pd}(z_0, \ldots, z_n, z_0', \ldots, z_n')$ holds.

The following lemma will be proved in the next section:

**Lemma 4** *Assume that $(K, v)$ is a valued field of characteristic $p > 0$ with $t \in K$ such that $vt$ is the smallest positive element in the value group $vK$. Then the sum of the groups $\wp(K), tK^p, \ldots, t^{p-1}K^p$ is pseudo direct. That is, $\mathrm{PD}(\wp(X), tX^p, \ldots, t^{p-1}X^p)$ holds in $(K, v)$.*

We need one further notion, which will play a key role in our results. A subset $S$ of a valued field $(K, v)$ will be called an **optimal approximation subset in** $(K, v)$ if for every $z \in K$ there is some $y \in S$ such that $v(z - y) = \max\{v(z - x) \mid x \in S\}$, i.e., if the following holds in $(K, v)$:

$$\forall Z \, \exists Y \in S \, \forall X \in S \ \mathcal{O}((Z - Y)/(Z - X)) . \tag{9}$$

(Recall that with our way of writing valuations, two points $x, y$ are the closer to each other, the bigger the value $v(x-y)$ is.) Note that (9) is an $\mathcal{L}(K)$-sentence if $S$ is $\mathcal{L}(K)$-definable.

For additive polynomials $f_0, \ldots, f_n \in K[X]$, we define:

$$\mathrm{OA}(f_0, \ldots, f_n) \ :\Leftrightarrow \ \text{the sum of the images of } f_0, \ldots, f_n$$
$$\text{is an optimal approximation subset.}$$

Since the subgroup $f_0(K) + \ldots + f_n(K)$ of $(K, +)$ is $\mathcal{L}(K)$-definable, $\mathrm{OA}(f_0, \ldots, f_n)$ is in fact an $\mathcal{L}(K)$-sentence. If $K$ is infinite (which we will always assume here, and which is automatic if $v$ is non-trivial), then also the fact that a polynomial $f$ is additive can be stated by an $\mathcal{L}(K)$-sentence:

$$\mathrm{ADD}(f) \ :\Leftrightarrow \ \forall X \forall Y \ f(X + Y) = f(X) + f(Y) .$$

Therefore, also the following is an $\mathcal{L}(K)$-sentence:

$$\left( \bigwedge_{i=0}^n \mathrm{ADD}(f_i) \wedge \mathrm{PD}(f_0, \ldots, f_n) \right) \Rightarrow \mathrm{OA}(f_0, \ldots, f_n) . \tag{10}$$

It asserts that if the given polynomials $f_0, \ldots, f_n$ are additive and the sum of their images is pseudo direct, then this sum is an optimal approximation subset.

The constants from $K$ can be removed by quantifying over the coefficients of the polynomials $f_0, \ldots, f_n$. By this method, for every $n \in \mathbb{N}$ we can get an elementary $\mathcal{L}$-sentence talking about at most $n + 1$ additive polynomials of degrees at most $p^n$. We obtain a recursive $\mathcal{L}$-axiom scheme expressing the following elementary property:

**(PDOA)**  for every $n \in \mathbb{N}$ and every choice of additive polynomials $f_0, \ldots, f_n$,
$$\mathrm{PD}(f_0, \ldots, f_n) \;\Rightarrow\; \mathrm{OA}(f_0, \ldots, f_n) \,.$$

One of our main results is:

**Theorem 5** *(PDOA) holds in every maximal field.*

We will give a proof in Section 2 below. (For maximal fields of characteristic 0, the theorem is trivial because then the only additive polynomials are of the form $cX$.)

Keeping some faith in our original sentence (5), let us observe:

**Lemma 6** *If $(K,v)$ satisfies axiom system (1) and (PDOA), then (5) holds in $(K,v)$ and $K$ satisfies (6) and (7).*

The proof will be given in the next section. Theorem 5 and Lemma 6 yield:

**Corollary 7** *If $(K,v)$ is a maximal field which satisfies axiom system (1), then (5) holds in $(K,v)$ and $K$ satisfies (6) and (7).*

Let us take advantage of the fact that (PDOA) is already formalized in $\mathcal{L}$, without needing the constant $t$. So far, we have kept secret the fact that we are much more interested in the $\mathcal{L}$-axiom system

$$\left.\begin{array}{l} (K,v) \text{ is a henselian defectless valued field} \\ K \text{ is of characteristic } p \\ vK \text{ is a } \mathbb{Z}\text{-group} \\ Kv = \mathbb{F}_p \end{array}\right\} \tag{11}$$

rather than in the $\mathcal{L}(t)$-axiom system (1). We only formulated (1) to show what the sentence (5) tells us
bout it. But now, we can derive:

**Theorem 8** *The axiom system (11) is not complete.*

Indeed, Theorem 5 shows that the model $(\mathbb{F}_p((t)), v_t)$ satisfies (PDOA), whereas Lemma 6 shows that the model $(L,v)$ given in Theorem 3 cannot satisfy (PDOA).

Now our main open question is:

**Is the axiom system (11) + (PDOA) complete?**

If this is the case, then it will also follow that $\mathrm{Th}(\mathbb{F}_p((t)), v_t)$ is decidable. We do not know an answer to this question. But we know that (PDOA) plays an important role in the structure theory of valued function fields. In fact, it admits to derive structure theorems of the same sort as we employed to prove the Ax–Kochen–Ershov principles for the elementary class of tame fields (cf. [K1], [K2], [K7]). Also, we can show that valued fields $(K,v)$ satisfying (11) + (PDOA) will satisfy the Ax-Kochen-Ershov principle with $\prec_\exists$ for arbitrary extensions $(L,v)$, provided that the extension $L|K$ is of transcendence degree 1. However, this needs an abundance of valuation theoretical machinery. The reason is that (PDOA) does not have as nice properties as "henselian" (or "tame"). Let us present one of the problems. It is a well known fact that a relatively algebraically closed subfield of a henselian field is again henselian. (The same holds for "tame" in the place of "henselian" if the extension is immediate.) But now consider an arbitrary maximal immediate extension $(M, v)$ of the field $(L, v)$ which is given in Theorem 3. By Theorem 5, (PDOA) holds in $(M, v)$. But it does not hold in $(L, v)$. On the other hand, the fact that $(L, v)$ is henselian defectless yields that $(L, v)$ is algebraically maximal. Therefore, it is relatively algebraically closed in $M$. Hence:

**Theorem 9** *There is an immediate extension $(L, v) \subset (M, v)$ of henselian defectless fields such that $L$ is relatively algebraically closed in $M$ and (PDOA) holds in $(M, v)$, but not in $(L, v)$.*

Another important property of "henselian" is: if $(K, v)$ is henselian, then so is each of its algebraic extensions. Also, "defectless" carries over to every finite extension (but not to every algebraic extension in general). So the following yet unanswered questions arise:

**Does (PDOA) carry over to finite extensions or even to algebraic extensions? What are the "algebraic properties" of (PDOA)?**

In the possible absence of uniqueness of maximal immediate extensions, e.g. for elementary classes containing $(\mathbb{F}_p((t)), v_t)$, one has to employ new ideas for the proof of Ax–Kochen–Ershov principles. Our proof for the case of tame fields profits from the fact that every extension of tame fields of finite transcendence degree can be split into an "anti-immediate" extension plus a tower of immediate extensions of tame fields of transcendence degree 1. In [K1] we proved a model theoretical result which takes care of the anti-immediate extension. For the immediate extensions of transcendence degree 1, we employ our structure theory for valued function fields. The reduction to transcendence degree 1 shows that in a certain sense, the model theoretical behaviour of tame fields (and of the other fields which we cited in (8)) is "one-dimensional". In contrast to this, (PDOA) tells us something about the correlation between several polynomials, and this is a "higher–dimensional information". Indeed, one can read off from Theorem 9 that for the case of $\mathbb{F}_p((t))$ a reduction to the case of transcendence degree 1 is much harder or even impossible (at least if one wants to remain in the given elementary class of valued fields). In fact, by a modification of our basic construction, we will show in Section 4 that for every $n \in \mathbb{N}$, we can construct $(L, v)$ in such a way that in addition to the assertions of Theorem 3, the following holds:

If $(L'|L, v)$ is an extension such that $L'v = Lv$ and (PDOA) holds in $(L', v)$, then $\mathrm{trdeg}\, L'|L \geq n$.

If we do not insist in $L'|L$ having finite transcendence degree, then we can even get that $\mathrm{trdeg}\, L'|L$ must be infinite.

For the conclusion of this section, let us think about three possible generalizations of Theorem 5:

1) It seems not unlikely to prove that already $\mathrm{OA}(f_0, \ldots, f_n)$ holds in every maximal field, for all additive polynomials $f_0, \ldots, f_n$. A possible way to prove this could be to show that for every choice of additive polynomials $f_0, \ldots, f_n$ there are additive polynomials $g_0, \ldots, g_m$ such that $\mathrm{PD}(g_0, \ldots, g_m)$ holds and $f_0(K) + \ldots + f_n(K) = g_0(K) + \ldots + g_m(K)$.

Let us call the group $f_0(K)+\ldots+f_n(K)$ a **polygroup**. So the generalization would state that every polygroup in a maximal field is an optimal approximation subset.

2) A polynomial $f(X_1,\ldots,X_n) \in K[X_1,\ldots,X_n]$ is called **additive** if it induces an additive map on $L^n$ for every extension field $L$ of $K$. With $f_i(X_i) = f(0,\ldots,0,X_i,0\ldots,0)$, it follows by additivity that $f(a_1,\ldots,a_n) = f_1(a_1)+\ldots+f_n(a_n)$ for all $(a_1,\ldots,a_n) \in L^n$. Since the polynomials $f_i$ are additive in one variable, we find that the polygroups in $K$ are precisely the images of the additive polynomials in several variables on $K$. Hence, the generalization indicated in 1) would actually state that the image of every additive polynomial in several variables on a maximal field is an optimal approximation subset.

3) Perhaps, the image of *every* polynomial in several variables on a maximal field is an optimal approximation subset. This would be an amazing generalization of Theorem 5 and of Lemma 12 of the next section. Our hope is that one could derive such a result from generalization 2) by approximating arbitrary polynomials by suitably chosen additive polynomials. For polynomials in one variable, something like this can be done by building on Kaplansky's work [KA].

## 2  Spherical completeness and optimal approximation

In the following, we will give the proof of Theorem 5. We need some further definitions. They can be given already in the context of ultrametric spaces, but here we will give them for subsets $S$ of valued fields $(K,v)$. A **closed ball in** $S$ is a set of the form $B_\gamma(a,S) = \{x \in S \mid v(a-x) \geq \gamma\}$ for $a \in S$ and $\gamma \in v(S-S) := \{v(s-s') \mid s,s' \in S\}$. A **nest of (closed) balls B** is a nonempty collection of closed balls such that each two balls in **B** have a nonempty intersection. By the ultrametric triangle law it follows that the balls in **B** are linearly ordered by inclusion. Now $(S,v)$ is called **spherically complete** if every nest of balls **B** in $S$ has a nonempty intersection: $\bigcap_{B \in \mathbf{B}} B \neq \emptyset$. It is easy to prove that $(S,v)$ is spherically complete if and only if every pseudo–convergent sequence in $(S,v)$ has a limit in $S$ (see [KA] or [K2] for these notions). Therefore, the following characterization of maximal fields is a direct consequence of Theorem 4 of [KA]:

**Theorem 10** *A valued field $(K,v)$ is maximal if and only if it is spherically complete.*

On the other hand, we have:

**Lemma 11** *Take any subset $S$ of the additive group of a valued field $(K,v)$. If $(S,v)$ is spherically complete, then $S$ is an optimal approximation subset in $(K,v)$.*

Proof:   Assume that $S$ is not an optimal approximation subset in $(K,v)$. Then there is an element $z \in K$ such that for every $y \in S$ there is some $x \in S$ satisfying that $v(z-x) > v(z-y)$. Note that by the ultrametric triangle law, the latter implies that $v(z-y) = v(x-y) \in v(S-S)$. From this and the fact that $S \cap B_{v(z-y)}(z,K) = B_{v(z-y)}(y,S)$, it follows that
$$\{\,B_{v(z-y)}(y,S) \mid y \in S\,\}$$
is a nest of balls in $(S,v)$. Take any $a \in S$ and choose $b \in S$ such that $v(z-b) > v(z-a)$. Then $a \notin B_{v(z-b)}(z,K)$. Hence the nest has an empty intersection, showing that $(S,v)$ is not spherically complete.  □

Now a natural question is: if $(K,v)$ is spherically complete and $f$ is an additive polynomial, does it follow that $(f(K),v)$ is spherically complete? In fact, this is true for *every* polynomial (the proof of the next two lemmas will be given in Section 5):

**Lemma 12** *If $(K,v)$ is spherically complete, then for every $f \in K[X]$, $(f(K),v)$ is spherically complete and therefore, $f(K)$ is an optimal approximation subset of $(K,v)$.*

Using Kaplansky's results together with the methods developed in Section 5), we can prove even more:

**Lemma 13** *If $(K,v)$ is algebraically maximal, then for every $f \in K[X]$, $f(K)$ is an optimal approximation subset of $(K,v)$.*

Now this exhibits an intriguing fact: if we have additive polynomials $f_0, \ldots, f_n$ on $K$ and $(K,v)$ is henselian defectless, then the $f_i(K)$ are optimal approximation subgroups, but their sum is not necessarily an optimal approximation subgroup, even if it is pseudo direct. By virtue of Lemma 15 below, the field $(L,v)$ of Theorem 3 with the additive polynomials $\wp(X), tX^p, \ldots, t^{p-1}X^p$ is an example for this. The situation changes when the subgroups are spherically complete:

**Theorem 14** *Let $G_0, \ldots, G_n$ be spherically complete subgroups of an arbitrary valued abelian group $(G,v)$. If their sum is pseudo direct, then it is spherically complete and hence an optimal approximation subset of $(G,v)$.*

The proof is given in [K4], using a theorem about maps on spherically complete ultrametric spaces. (It seems unlikely that the theorem works without any condition on the sum of the $G_i$'s. But we do not know of any counterexample.)

Now Theorem 5 follows from Theorem 10, Lemma 12 and Theorem 14. We note that all this works as well for arbitrary definable additive maps in the place of additive polynomials. However, we do not know of any such maps which would provide subgroups essentially different from polygroups. More generally, we ask:

**Do there exist definable subgroups in valued fields of positive characteristic which are essentially different from polygroups? Are they optimal approximation subgroups? Do they carry other (independent) valuation theoretical properties? Are there polygroups which are not representable as pseudo direct sums but are optimal approximation subgroups?**

Let us note that it makes no essential difference to add the "trivial" subgroups like $\mathcal{O}$, $\mathcal{M}$ or other balls around 0. Optimal approximation assertions about groups obtained in such a way from polygroups are consequences of (PDOA).

## 3 Valuation independence and pseudo direct sums

Take any valued field extension $(K|K',v)$. The elements $c_0, \ldots, c_m \in K \setminus \{0\}$ will be called $K'$-**valuation independent** if for every choice of elements $d_0, \ldots, d_m \in K'$, the following holds:
$$v(c_0 d_0 + \ldots + c_m d_m) = \min_{0 \leq i \leq m} v c_i d_i .$$

In particular, if $d_k \neq 0$ for at least one $k$, then $c_k d_k \neq 0$ and thus, $v(c_0 d_0 + \ldots + c_m d_m) \leq v c_k d_k < \infty$ which shows that $c_0 d_0 + \ldots + c_m d_m \neq 0$. Hence if $c_0, \ldots, c_m$ are $K'$-valuation independent, then they are $K'$-linearly independent.

**Lemma 15** *Take a valued field $(K, v)$ of characteristic $p > 0$ with $K^p$-valuation independent elements $c_0, \ldots, c_m$, where $c_0 = 1$. Then $\mathrm{PD}(\wp(X), c_1 X^p, \ldots, c_m X^p)$ holds in $(K, v)$.*

Proof:   We set $f_0(X) = \wp(X)$ and $f_i(X) = c_i X^p$ for $1 \leq i \leq m$. For $x_0, \ldots, x_m \in K$,
$$f_0(x_0) + \ldots + f_m(x_m) = -x_0 + c_0 x_0^p + \ldots + c_m x_m^p \,.$$

If $v x_0 > v(c_0 x_0^p + \ldots + c_m x_m^p)$, then
$$\begin{aligned} v(f_0(x_0) + \ldots + f_m(x_m)) &= \min\{v x_0, v(c_0 x_0^p + \ldots + c_m x_m^p)\} \\ &= v(c_0 x_0^p + \ldots + c_m x_m^p) = \min_{0 \leq i \leq m} v c_i x_i^p = \min_{0 \leq i \leq m} v f_i(x_i) \,,\end{aligned}$$

where the last equality holds since $v x_0 > v c_0 x_0^p$ implies that $v c_0 x_0^p = v f_0(x_0)$.

Now assume that $v \sum_{i=0}^{m} f_i(x_i) > \min_i v f_i(x_i)$. Then by what we just have shown,
$$v x_0 \leq v(c_0 x_0^p + \ldots + c_m x_m^p) = \min_{0 \leq i \leq m} v c_i x_i^p \leq v c_0 x_0^p = p v x_0 \,.$$

But $v x_0 \leq p v x_0$ can only hold if $v x_0 \geq 0$, in which case also $v f_0(x_0) \geq 0$. We also find that $0 \leq v x_0 \leq \min_i v c_i x_i^p \leq v c_j x_j^p = v f_j(x_j)$ for all $j \geq 1$. Hence, $\min_i v f_i(x_i) \geq 0$. Now it follows from our assumption that $v \sum_{i=0}^{m} f_i(x_i) > 0$. We set $y_0 = -\sum_{i=0}^{m} f_i(x_i)$. Observe that $v y_0 > 0$ implies that $v y_0^p > v y_0$. Hence,
$$v\left(\sum_{i=0}^{m} f_i(x_i) - \wp(y_0)\right) = v y_0^p > v y_0 = v \sum_{i=0}^{m} f_i(x_i) \,.$$

Taking $y_i = 0$ for $i \geq 1$, we obtain that
$$\mathrm{pd}(f_0(x_0), \ldots, f_m(x_m), f_0(y_0), \ldots, f_m(y_m))$$

holds. □

If in the situation of this lemma, $(K, v)$ is henselian, then we can even get that $\sum_{i=0}^{m} f_i(y_i) = \sum_{i=0}^{m} f_i(x_i)$. Indeed, using that $\mathcal{M} \subset \wp(K)$, in the second part of the proof we just have to choose $y_0 \in K$ such that $\wp(y_0) = \sum_{i=0}^{m} f_i(x_i)$.

**Lemma 16** *Assume that $(K, v)$ is a valued field of characteristic $p > 0$ with $t \in K$ such that $vt$ is the smallest positive element in the value group $vK$. Then the elements $1, t, t^2, \ldots, t^{p-1}$ are $K^p$-valuation independent.*

Proof: For every choice of elements $d_0, \ldots, d_{p-1}$ we have that $vt^i d_i^p \in ivt + pvK$. As $vt$ is the smallest positive element of $vK$ by assumption, the cosets $pvK, vt + pvK, 2vt + pvK, \ldots, (p-1)vt + pvK$ are all distinct. This shows that $vt^i d_i^p \neq vt^j d_j^p$ for $0 \leq i < j \leq p-1$. Hence, $v(d_0 + td_1 + \ldots + t^{p-1} d_{p-1}) = \min_{0 \leq i \leq p-1} vt^i d_i$. □

Now Lemma 4 follows from Lemmas 15 and 16.

A valued field $(K, v)$ is called **inseparably defectless** if the fundamental equality holds for every finite purely inseparable extension. We will need the following characterization of inseparably defectless fields, which was proved by F. Delon [D] (see also [K2]):

**Lemma 17** *Take a valued field $(K, v)$ of characteristic $p > 0$ such that $(vK : pvK) < \infty$ and $[Kv : (Kv)^p] < \infty$. Then $(K, v)$ is inseparably defectless if and only if*

$$[K : K^p] = (vK : pvK)[Kv : (Kv)^p] \,. \tag{12}$$

From this lemma we obtain:

**Lemma 18** *Assume that $(K, v)$ is an inseparably defectless valued field of characteristic $p > 0$ with $t \in K$ such that $vt$ is the smallest positive element in the value group $vK$. Assume further that $vK$ is a $\mathbb{Z}$-group and that $Kv$ is perfect. Then $1, t, t^2, \ldots, t^{p-1}$ is a basis of $K|K^p$.*

Proof: By Lemma 16 and our remark in the beginning of this section we know that $1, t, t^2, \ldots, t^{p-1}$ are $K^p$-linearly independent. By the foregoing lemma, (12) holds. Since $vK$ is a $\mathbb{Z}$-group, we have that $(vK : pvK) = p$. Since $Kv$ is perfect, we have that $[Kv : (Kv)^p] = 1$. Thus, $[K : K^p] = p$, which shows that $1, t, t^2, \ldots, t^{p-1}$ is a basis of $K|K^p$. □

**Lemma 19** *Let the assumptions be as in Lemma 18. Take $z \in K$ and assume that the set*

$$\{v(z - y) \mid y \in \wp(K) + tK^p + \ldots + t^{p-1}K^p\} \tag{13}$$

*admits a maximum. Then this maximum is either $0$ or $\infty$ (the latter meaning that $z$ lies in $\wp(K) + tK^p + \ldots + t^{p-1}K^p$).*

Proof: Assume that $y_0 \in K$ is such that $v(z - y_0)$ is the maximum of (13). After replacing $z$ by $z - y_0$ we can assume that $y_0 = 0$.

Suppose that $vz > 0$. Then $vz = \infty$ since otherwise, we could set

$$y := -z^p + z = (-z)^p - (-z) + t \cdot 0 + \ldots + t^{p-1} \cdot 0 \in \wp(K) + tK^p + \ldots + t^{p-1}K^p$$

to obtain that $v(z - y) = vz^p > vz$, a contradiction.

Now suppose that $vz < 0$. We have to deduce a contradiction from this assumption. By Lemma 18, we can write

$$z = b_0^p + tb_1^p + \ldots + t^{p-1}b_{p-1}^p \quad \text{with } b_0, \ldots, b_{p-1} \in K \,.$$

By Lemma 16 we have that
$$\min_{0\leq i\leq p-1} vt^i b_i^p \;=\; v(b_0^p + tb_1^p + \ldots + t^{p-1}b_{p-1}^p) \;=\; vz \;<\; 0 \;.$$

Hence if $vb_0 < 0$, then
$$vb_0 \;>\; pvb_0 \;=\; vb_0^p \;\geq\; \min_{0\leq i\leq p-1} vt^i b_i^p \;=\; vz \;.$$

On the other hand, if $vb_0 \geq 0$, then $vb_0 > vz$, too. Hence in every case,
$$v(z - (\wp(b_0) + tb_1^p + \ldots + t^{p-1}b_{p-1}^p)) \;=\; vb_0 \;>\; vz \;,$$

a contradiction to the maximality of $vz$. $\square$

**Proof of Lemma 6**: Assume that $(K, v)$ satisfies axiom system (1) and (PDOA). By Lemma 15 in connection with Lemma 18, $\mathrm{PD}(\wp(X), tX^p, \ldots, t^{p-1}X^p)$ holds in $(K, v)$. Thus, (PDOA) yields that $\mathrm{OA}(\wp(X), tX^p, \ldots, t^{p-1}X^p)$ holds in $(K, v)$. Take $z \in K$ and suppose that $z \notin \wp(K) + tK^p + \ldots + t^{p-1}K^p$. Then by Lemma 19, there is some $y \in \wp(K) + tK^p + \ldots + t^{p-1}K^p$ such that $v(z - y) = 0$. Since $Kv = \mathbb{F}_p$, there is some $j \in \mathbb{F}_p$ such that $v(z - y - j) > 0$. Again by Lemma 19, we obtain that $z - y - j \in \wp(K) + tK^p + \ldots + t^{p-1}K^p$. This proves that $K$ satisfies (7). Hence, it also satisfies (6), and (5) holds in $(K, v)$. $\square$

## 4 An example and its consequences

We need some preparations. The **rank** of $(K, v)$ is the number of proper convex subgroups of the value group $vK$ (if finite); $(K, v)$ has rank 1 if and only if $vK$ is archimedean, i.e., embeddable in the additive group of the reals. If $(K, v)$ has rank $n$, then $v$ is the composition of $n$ valuations of rank 1.

Here is a well-known fact about pseudo–convergent sequences. Unfortunately, it is not explicitly stated in [KA].

**Lemma 20** *Assume that $(a_\nu)_{\nu<\lambda}$ is a pseudo–convergent sequence in a valued field $(L, v)$ (where $\lambda$ is a limit ordinal). If $b$ is not a limit of this sequence, then there is some $\nu_0 < \lambda$ such that for all $\nu \geq \nu_0$, $\nu < \lambda$,*
$$v(b - a_\nu) \;<\; v(a_{\nu_0+1} - a_{\nu_0}) \;.$$

Proof: By definition, $b$ is not a limit of $(a_\nu)_{\nu<\lambda}$ if and only if $v(b - a_\mu) \neq v(a_{\mu+1} - a_\mu)$ for some $\mu < \lambda$. We set $\nu_0 := \mu + 2$. Take any $\nu \geq \nu_0$ such that $\nu < \lambda$. Then by Lemma 2 of [KA] and by the definition of pseudo–convergent sequences,
$$v(a_\nu - a_{\mu+1}) \;=\; v(a_{\mu+2} - a_{\mu+1}) \;>\; v(a_{\mu+1} - a_\mu) \;.$$

So we obtain that
$$\begin{aligned} v(b - a_\nu) \;&=\; \min\{v(b - a_\mu), v(a_\mu - a_{\mu+1}), v(a_{\mu+1} - a_\nu)\} \\ &=\; \min\{v(b - a_\mu), v(a_\mu - a_{\mu+1})\} \;\leq\; v(a_\mu - a_{\mu+1}) \\ &<\; v(a_{\mu+2} - a_{\mu+1}) \;=\; v(a_{\nu_0+1} - a_{\nu_0}) \;, \end{aligned}$$

where the first equality holds since the values $v(b - a_\mu)$, $v(a_\mu - a_{\mu+1})$ and $v(a_{\mu+1} - a_\nu)$ are distinct. □

We will need the following characterization of henselian defectless fields; for a proof, see [K2].

**Theorem 21** *Let $(K,v)$ be an inseparably defectless field of characteristic $p > 0$. If in addition $(K,v)$ is algebraically maximal, then $(K,v)$ is henselian defectless.*

Now we are ready for the construction of a basic example, which we will then use to prove Theorem 3. Let $K$ be a field of characteristic $p > 0$. Further, assume that $m := [K : K^p]$ is finite, and choose a basis $c_0, \ldots, c_m$ of $K|K^p$ with $c_0 = 1$. We work in the power series field $K((s^{\mathbb{Q}}))$ with its canonical ($s$-adic) valuation $v_s$. As this is henselian, it contains the henselization of the subfield $K(s^{1/n} \mid n \in \mathbb{N}, (p,n) = 1)$ with respect to (the restriction of) $v_s$. We will denote this henselization by $L_1$. We have that

$$v_s L_1 = \sum_{n \in \mathbb{N},\, (p,n)=1} \frac{1}{n} \mathbb{Z} . \tag{14}$$

In particular, $1/q \in v_s L_1$ and $s^{1/q} \in L_1$ for every prime number $q \neq p$.

We take an ascending sequence of prime numbers $q_j$, $j \in \mathbb{N}$, such that

$$p^{j+1} < q_j \quad \text{for all } j \in \mathbb{N} . \tag{15}$$

In particular, $p < q_j$ and thus $s^{1/q_j} \in L_1$ for all $j \in \mathbb{N}$. Now assume that $\zeta$ is a limit of the pseudo–convergent sequence

$$\left( \sum_{j=1}^{k} s^{-1/q_j} \right)_{k \in \mathbb{N}} \tag{16}$$

in some extension of $(L_1, v_s)$. Using the method employed in Example 16.1 of [K5] one shows by use of Hensel's Lemma that

$$v_s K(s, \zeta) = \sum_{j \in \mathbb{N}} \frac{1}{q_j} \mathbb{Z} .$$

Condition (15) yields that the sequence $(q_j)$ contains infinitely many primes; consequently, $(v_s K(s, \zeta) : v_s K(s)) = (v_s K(s, \zeta) : \mathbb{Z})$ is not finite. By virtue of the fundamental inequality, this shows that $\zeta$ must be transcendental over $K(s)$ and thus also over its algebraic extension $L_1$. By virtue of Theorem 3 of [KA], this proves that the pseudo–convergent sequence $(\sum_{j=1}^{k} s^{-1/q_j})_{k \in \mathbb{N}}$ in $(L_1, v_s)$ cannot be of algebraic type; hence it must be of transcendental type.

We will now construct a purely inseparable algebraic extension $L_2$ of $L_1$ such that $c_0, \ldots, c_m$ is again a basis of $L_2|L_2^p$. We define recursively

$$\xi_1 = s^{-1/p} \quad \text{and} \quad \xi_{j+1} = (\xi_j - c_1 s^{-p/q_j})^{1/p} . \tag{17}$$

Since $v_s$ is trivial on $K$, we have that $v_s c_1 = 0$. Using this and (15), one shows by induction on $j$ that

$$v_s \xi_j \;=\; -\frac{1}{p^j} \;<\; -\frac{p}{q_j} \;=\; v_s(c_1 s^{-p/q_j}) \;<\; 0 \quad \text{for all } j \in \mathbb{N} \,. \tag{18}$$

We put
$$L_2 \;:=\; L_1(\xi_j \mid j \in \mathbb{N}) \,.$$

To prove that $c_0, \ldots, c_m$ is a $p$-basis of $L_2$, take $a \in L_2$. Then $a \in L_1(\xi_1, \ldots, \xi_k) = L_1(\xi_k)$ for a suitable $k \in \mathbb{N}$. Now one deduces by induction that $c_\mu \xi_k^\nu$, $0 \leq \mu \leq m$, $0 \leq \nu < p$, is a $p$-basis for $L_1(\xi_k)$ and that

$$\xi_k \;=\; \xi_{k+1}^p + c_1 s^{-p/q_j} \in L_1(\xi_{k+1})^p + c_1 L_1(\xi_{k+1})^p \;\subset\; K . L_2^p \,.$$

This shows that

$$a \in \sum_{\mu,\nu} c_\mu \xi_k^\nu L_1(\xi_k)^p \subset K.L_2^p \;=\; c_0 L_2^p + c_1 L_2^p + \ldots + c_m L_2^p \,.$$

Hence, $1, c_1, \ldots, c_m$ is a $p$-basis of $L_2$.

Since the extension $L_2 | L_1$ is purely inseparable, there is a unique extension $w$ of $v_s$ to $L_2$. (Note that we are now working outside of $K((s^\mathbb{Q}))$). For each $j$ we have that $\frac{1}{p^j} = w\xi_j \in wL_1(\xi_j)$ and thus, $(wL_1(\xi_j) : v_s L_1) \geq p^j = [L_1(\xi_j) : L_1]$. By the fundamental inequality, $[L_1(\xi_j) : L_1] \geq (wL_1(\xi_j) : v_s L_1)$. Hence, $[L_1(\xi_j) : L_1] = (wL_1(\xi_j) : v_s L_1)$, and

$$wL_1(\xi_j) \;=\; v_s L_1 + \frac{1}{p^j}\mathbb{Z} \,. \tag{19}$$

Again by the fundamental inequality it follows that $L_1(\xi_j)w = L_1 v_s = K$ and therefore,

$$L_2 w \;=\; \bigcup_{j \in \mathbb{N}} L_1(\xi_j) w \;=\; K \,.$$

By (19) and (14),

$$wL_2 \;=\; \bigcup_{j \in \mathbb{N}} wL_1(\xi_j) \;=\; \bigcup_{j \in \mathbb{N}} \left( v_s L_1 + \frac{1}{p^j}\mathbb{Z} \right) \;=\; \mathbb{Q} \,. \tag{20}$$

Now we choose $(L, w)$ to be a maximal immediate algebraic extension of $(L_2, w)$ (which exists by Zorn's Lemma since its cardinality is bounded by that of the algebraic closure of $L_2$). Then

$$\begin{aligned}
[L : L^p] \;&\leq\; [L_2 : L_2^p] \;=\; m + 1 \;=\; [K : K^p] \;=\; [L_2 w : (L_2 w)^p] \\
&=\; [L_2 w : (L_2 w)^p] \cdot (\mathbb{Q} : p\mathbb{Q}) \;=\; [L_2 w : (L_2 w)^p] \cdot (wL_2 : pwL_2) \\
&=\; [Lw : (Lw)^p] \cdot (wL : pwL) \;\leq\; [L : L^p] \,.
\end{aligned}$$

Hence, equality holds everywhere. By Lemma 17, the last equality implies that $(L, w)$ is inseparably defectless. Since $(L, w)$ is a maximal immediate algebraic extension and thus algebraically maximal, Theorem 21 shows that $(L, w)$ is a henselian defectless field.

We set
$$x := s^{-1} .$$

Assume that there is an extension $(L'|L, v)$ such that $L'v = Lv$, and that there exist elements $x_0, x_1, \ldots, x_m, y \in L'$ such that

$$x = y + x_0^p - x_0 + c_1 x_1^p + \ldots + c_m x_m^p \quad \text{with } wy \geq 0 . \tag{21}$$

We wish to show that then $x_1$ must be a limit of the pseudo–convergent sequence (16), which yields that $x_1$ is transcendental over $L$. This in turn shows that (21) cannot hold in $L$.

Suppose that $x_1$ is not a limit of (16). Then by Lemma 20, there exists some $k_0 \in \mathbb{N}$ such that for all $k \geq k_0$,

$$w\left(x_1 - \sum_{j=1}^{k} s^{-1/q_j}\right) < vs^{-1/q_{k_0+1}} = -\frac{1}{q_{k_0+1}} . \tag{22}$$

We can choose $k$ as large as to also guarantee that $p^k > q_{k_0+1}$, that is,

$$-\frac{1}{q_{k_0+1}} < -\frac{1}{p^k} = w\xi_k . \tag{23}$$

We set

$$\tilde{x}_0 := x_0 - \sum_{j=1}^{k} \xi_j \quad \text{and} \quad \tilde{x}_1 := x_1 - \sum_{j=1}^{k-1} s^{-1/q_j} . \tag{24}$$

According to (22) and (23), we have that

$$pw\tilde{x}_1 < w\tilde{x}_1 < w\xi_k < 0 . \tag{25}$$

Now we compute

$$\begin{aligned}
\tilde{x}_0^p - \tilde{x}_0 &= x_0^p - x_0 + \left(-\sum_{j=1}^{k} \xi_j\right)^p + \sum_{j=1}^{k} \xi_j \\
&= x_0^p - x_0 - \xi_1^p - \sum_{j=1}^{k-1} (\xi_{j+1}^p - \xi_j) + \xi_k \\
&= x_0^p - x_0 - x + \sum_{j=1}^{k-1} c_1 s^{-p/q_j} + \xi_k \\
&= \xi_k - y - (c_1 \tilde{x}_1^p + c_2 x_2^p + \ldots + c_m x_m^p) .
\end{aligned}$$

Since $w\xi_{k+1} < 0$ and $wy \geq 0$, and by virtue of (25), we have that

$$0 > w(\xi_k - y) = w\xi_k > pw\tilde{x}_1 = w\tilde{x}_1^p \geq \min\{w\tilde{x}_1^p, wx_2^p, \ldots, wx_m^p\} =: \alpha .$$

We set $\tilde{x}_i := x_i$ for $2 \leq i \leq m$ and take $i_1, \ldots, i_\ell \in \{1, \ldots, m\}$ to be all indices $i$ for which $w\tilde{x}_i^p = \alpha$. Then $w(\tilde{x}_{i_\nu}^p/\tilde{x}_{i_1}^p) = 0$ and $w((\xi_{k+1} - y)/\tilde{x}_{i_1}^p) > 0$. Therefore, and since the elements $1, c_1, \ldots, c_m \in K$ are linearly independent over $K^p = (Lw)^p = (L'w)^p$,

$$\frac{\tilde{x}_0^p - \tilde{x}_0}{\tilde{x}_{i_1}^p} w = \left(c_{i_1} \frac{\tilde{x}_{i_1}^p}{\tilde{x}_{i_1}^p} + \ldots + c_{i_\ell} \frac{\tilde{x}_{i_\ell}^p}{\tilde{x}_{i_1}^p}\right) w = c_{i_1} + c_{i_2}\left(\frac{\tilde{x}_{i_2}}{\tilde{x}_{i_1}} w\right)^p + \ldots + c_{i_\ell}\left(\frac{\tilde{x}_{i_\ell}}{\tilde{x}_{i_1}} w\right)^p \notin (L'w)^p .$$

In particular, the residue is nonzero, which implies that
$$w(\tilde{x}_0{}^p - \tilde{x}_0) = \tilde{x}_{i_1}^p = \alpha < 0 \ .$$
This yields that $w\tilde{x}_0 < 0$. Consequently, $w\tilde{x}_0 > w\tilde{x}_0{}^p$ and thus,
$$\left(\frac{\tilde{x}_0}{\tilde{x}_{i_1}}w\right)^p = \frac{\tilde{x}_0{}^p}{\tilde{x}_{i_1}^p}w = \frac{\tilde{x}_0{}^p - \tilde{x}_0}{\tilde{x}_{i_1}^p}w \notin (L'w)^p \ .$$
This contradiction proves that $x_1$ must be a limit of (16).

Now we take $K$ to be any field of characteristic $p$ containing an element $t$ such that (2) holds (for example, we may take $K = \mathbb{F}_p(t)$, $K = \mathbb{F}_p(t)^h$ or $K = \mathbb{F}_p((t))$.) Then we can set $m = p - 1$ and $c_i = t^i$ for $0 \leq i \leq m$. We obtain that the existential sentence
$$\exists Y \exists X_0 \ldots \exists X_{p-1} \ x = Y + X_0^p - X_0 + tX_1^p + \ldots + t^{p-1}X_{p-1}^p \wedge \mathcal{O}(Y) \qquad (26)$$
does not hold in $(L, w)$. So we have proved:

**Theorem 22** *Let $K$ be any field of characteristic $p > 0$ containing an element $t$ such that $K = K^p \oplus tK^p \oplus \ldots \oplus t^{p-1}K^p$. Then there exists a henselian defectless field $(L, w)$, not satisfying property (5), of transcendence degree 1 over its embedded residue field $K$, having value group $wL = \mathbb{Q}$, and such that*
$$L = K.L^p = L^p \oplus tL^p \oplus \ldots \oplus t^{p-1}L^p \ .$$

Now we can give the

**Proof of Theorem 3:** We take $K$ to be the field $\mathbb{F}_p(t)^h$ or $\mathbb{F}_p((t))$, with $v_t$ the $t$-adic valuation on $K$. We denote by $v$ the composition $w \circ v_t$ of $w$ with $v_t$ on $L$; this can actually be viewed as an extension of $v_t$ to $L$. We note that $v$ is finer than $w$, that is, $\mathcal{O}_v \subset \mathcal{O}_w$. This means that $vy \geq 0$ implies $wy \geq 0$; therefore, since (5) doesn't hold for $(L, w)$, it doesn't hold for $(L, v)$.

We have mentioned already that both $(\mathbb{F}_p(t)^h, v_t)$ and $(\mathbb{F}_p((t)), v_t)$ are defectless fields. On the other hand, we know from our construction that $(L, w)$ is a henselian defectless field. Since the composition of henselian defectless valuations is again henselian defectless (cf. [K2]), it follows that for both choices of $K$, $(L, v)$ is a henselian defectless field. Since $wL = \mathbb{Q}$ and $v_t(Lw) = v_tK = \mathbb{Z}$, we have that $vt = v_tt$ is the smallest positive element of $vL$, $\mathbb{Z}vt$ is a convex subgroup of $vL$, and $vL/\mathbb{Z}vt \simeq \mathbb{Q}$. Hence, $vL$ is a Z-group. Further, $Lv = (Lw)v_t = Kv_t = \mathbb{F}_p$. By construction, $1, t, t^2, \ldots, t^{p-1}$ is a basis of $L|L^p$.

Finally, it remains to show that $L|K$ is regular. Take any finite extension $K'|K$ and an extension of $v_t$ from $L$ to $K'.L$. Since $(K, v_t)$ is henselian, the restriction of $v$ from $K'.L$ to $K'$ is the unique extension of $v_t$ from $K$ to $K'$. We set e := $(vK' : vK)$ and f := $[K'v : Kv_t]$. Since $(K, v_t)$ is defectless, we have that $[K' : K] = $ ef. As $v_tK = \mathbb{Z}v_tt$, there is some $t' \in K'$ such that $evt' = vt$; therefore, $t' \in K'.L$ yields that $(vK'.L : vL) \geq$ e. Since $Kv_t = \mathbb{F}_p = Lv$, we also find that $[(K'.L)v : Lv] = [(K'.L)v : \mathbb{F}_p] \geq [K'v : \mathbb{F}_p] = $ f. Thus,
$$[K' : K] = \text{ef} \leq (vK'.L : vL)[(K'.L)v : Lv] \leq [K'.L : L] \leq [K' : K] \ .$$

Therefore, equality must hold everywhere, showing that $L|K$ is linearly disjoint from $K'|K$. Since $K'|K$ was an arbitrary finite extension, this proves that $L|K$ is regular. □

**Remark 23** These examples also show that a field which is relatively algebraically closed in a henselian defectless field that satisfies (5) does itself not necessarily satisfy (5), even if the extension is immediate. Indeed, every maximal immediate extension of our examples $(L, v)$ or $(L, w)$ is a maximal field and thus satisfies (5) according to Theorem 7, and $L$ is relatively algebraically closed in every such extension since $(L, v)$ and $(L, w)$ are henselian defectless and thus algebraically maximal.

With the examples that we have constructed, we can even show a sharper result. Beforehand, we need two auxiliary lemmas. Note that it is easy to show that a pseudo–convergent sequence of transcendental type in $(k, v)$ will never have a limit in $k$.

**Lemma 24** *Take any henselian field $(k, v)$ and an immediate extension $(k(x)|k, v)$ such that $x$ is the limit of a pseudo–convergent sequence of transcendental type in $(k, v)$. Then $(k, v)$ is existentially closed in the henselization $(k(x), v)^h$ of $(k(x), v)$.*

Proof: Let $x$ be the limit of the pseudo–convergent sequence $(x_\nu)_{\nu<\lambda}$ of transcendental type in $(k, v)$. We take $(k^*, v^*)$ to be a $|k|^+$-saturated elementary extension of $(k, v)$. Every finite subset of the set $\{\mathcal{O}((X - a_\nu)/(a_{\nu+1} - a_\nu)) \mid \nu < \lambda\}$ of atomic $\mathcal{L}(k)$-sentences is satisfied by $a_\mu$ if $\mu$ is bigger than all $\nu$ appearing in that subset. Thus, there is an element $x^* \in k^*$ which simultaneously satisfies all sentences in this set. Then $x^*$ is a limit of $(x_\nu)_{\nu<\lambda}$ (we leave the easy proof to the reader). By Theorem 2 of [KA], $x \mapsto x^*$ induces a valuation preserving isomorphism of $k(x)$ onto $k(x^*)$. Since $(k^*, v^*)$ is an elementary extension of $(k, v)$ and "henselian" is an elementary property, also $(k^*, v^*)$ is henselian. Hence by the universal property of henselizations (cf. [R] or [K2]), this isomorphism can be extended to an embedding of $(k(x), v)^h$ in $(k^*, v^*)$. Now every existential $\mathcal{L}(k)$-sentence holding in $(k(x), v)^h$ carries over to $(k^*, v^*)$ through the embedding. Since $(k, v) \prec (k^*, v^*)$, it will therefore also hold in $(k, v)$. □

**Lemma 25** *Take a henselian field $(k, v)$, a polynomial $f \in k[X]$ of degree $p = \mathrm{char}\, kv$, and a root $a$ of $f$. Suppose that $(a_\nu)_{\nu<\lambda}$ is a pseudo–convergent sequence which does not fix the value of $f$ and has no limit in $(k, v)$. Then there is an immediate extension of $v$ from $k$ to $k(a)$ such that $a$ is a limit of $(a_\nu)_{\nu<\lambda}$ in $(k(a), v)$.*

Proof: We pick a polynomial $g \in k[X]$ of minimal degree with the property that $(a_\nu)_{\nu<\lambda}$ does not fix the value of $g$. Take a root $b$ of $g$. Then by Theorem 3 of [KA] there is an immediate extension of $v$ from $k$ to $k(b)$. Since $(k, v)$ is assumed to be henselian, we have e = f = g = 1. By the Lemma of Ostrowski (cf. [R] or [K2]), it follows that $\deg g = [k(b) : k] = [k(b) : k]/\mathrm{efg}$ is a power of $p$. This proves that $f$ is of minimal degree with the property that $(a_\nu)_{\nu<\lambda}$ does not fix the value of $f$. Hence, our assertion follows by a second application of Theorem 3 of [KA]. □

**Theorem 26** *Let $L$ be the field given by our construction. Then there exists a regular function field $F$ of transcendence degree 1 and generated by two elements over $L$ such that $L$ is not existentially closed in $F$ (in the language of rings), but $w$ and $v$ have immediate extensions from $L$ to $F$.*

Proof: We will show that the existential sentence

$$\exists Y \exists X_0 \ldots \exists X_{p-1} \quad x = X_0^p - X_0 + tX_1^p + \ldots + t^{p-1}X_{p-1}^p \tag{27}$$

holds already in $L(x_0, x_1)$, where

$$(L(x_0, x_1)|L, w)$$

is a regular immediate function field with $L(x_0, x_1)|L(x_1)$ an Artin-Schreier extension. As $L(x_0, x_1)w = Lw = K$, we can again define $v = w \circ v_t$ on $L(x_0, x_1)$. Then also $(L(x_0, x_1)|L, v)$ is immediate. This will imply the assertion of our lemma.

We take $x_1$ to be a transcendental element over $L$. Using Theorem 2 of [KA], we extend $w$ to $L(x_1)$ in such a way that $x_1$ becomes a limit of the pseudo–convergent sequence (16) and $(L(x_1)|L, w)$ is an immediate extension. Then we define a pseudo–convergent sequence $(a_k)_{k \in \mathbb{N}}$ by setting

$$a_k := \sum_{j=1}^{k} \xi_j . \tag{28}$$

Now we compute for all $k \in \mathbb{N}$, using (17) and (18):

$$w(a_k^p - a_k - (x - tx_1^p)) =$$
$$= w\left((\sum_{j=1}^{k} \xi_j)^p - \sum_{j=1}^{k} \xi_j - (x - tx_1^p)\right) = w\left(\xi_1^p + \sum_{j=1}^{k-1}(\xi_{j+1}^p - \xi_j) - \xi_k - (x - tx_1^p)\right)$$
$$= w\left(x - t\sum_{j=1}^{k-1} s^{-p/q_j} - \xi_k - (x - tx_1^p)\right) = w\left(t(x_1 - \sum_{j=1}^{k-1} s^{-1/q_j})^p - \xi_k\right)$$
$$= \min\left(wt(x_1 - \sum_{j=1}^{k-1} s^{-1/q_j})^p, w\xi_k\right) = \min\left(wts^{-p/q_k}, w\xi_k\right) = w\xi_k = -\frac{1}{p^k}$$

(where the first equality of the last line holds since $wts^{-p/q_k} \neq w\xi_k$). This shows that the pseudo–convergent sequence $(a_k)_{k \in \mathbb{N}}$ does not fix the value of the Artin–Schreier polynomial

$$X^p - X - (x - tx_1^p) . \tag{29}$$

Also, we see that a limit $x_0$ of $(a_k)_{k \in \mathbb{N}}$ in an arbitrary extension of $(L, v)$ will satisfy

$$w(x_0^p - x_0 - (x - tx_1^p)) > -\frac{1}{p^k} \quad \text{for all } k \in \mathbb{N} ,$$

whence

$$w(x_0^p - x_0 - (x - tx_1^p)) \geq 0 .$$

This means that the existential sentence (26) holds in $(L(x_0, x_1), w)$. From Lemma 24 we know that $(L, w)$ is existentially closed in the henselization $(L(x_1), w)^h$ the immediate rational function field $(L(x_1), w)$. So if $x_0$ were an element of this henselization, (26) would also hold in $(L, w)$, contrary to what we have already proved. This contradiction shows that the pseudo–convergent sequence $(a_k)_{k \in \mathbb{N}}$ has no limit in $(L(x_1), w)^h$. Hence by virtue of Lemma 25, if $x_0$ is any root of the polynomial (29), then there is an immediate extension of the valuation $w$ from $L(x_1)^h$ to $L(x_1)^h(x_0)$. Its restriction is an immediate extension of $w$ from $L(x_1)$ to the Artin-Schreier
xtension $L(x_0, x_1)$. Now (27) is satisfied in $L(x_0, x_1)$, as desired.

Using that $(F|L, v)$ is immediate, the regularity of $F|L$ can be shown in a similar way as the regularity of $L|K$ was shown in the proof of Theorem 3. This completes the proof of our lemma. □

By taking maximal immediate extensions of $(F, w)$ and of $(F, v)$, we obtain:

**Corollary 27** *There are maximal immediate extensions of $(L, w)$ and of $(L, v)$ in which $L$ is not existentially closed (already in the language of rings).*

Let us draw a further conclusion which is important for an application to the question of local uniformization (cf. [K3], [K5], [K8]). A valued function field $(F|L, v)$ of transcendence degree 1 is called **henselian rational** if there is some $x \in F^h$ such that $F^h = L(x)^h$. In [K1] we have proved that every immediate function field $(F|L, v)$ of transcendence degree 1 over a tame field $(L, v)$ is henselian rational (see also [K2] and [K7]). In contrast to this result, we have:

**Corollary 28** *The immediate function fields $(F|L, w)$ and $(F|L, v)$ of the foregoing theorem are not henselian rational.*

Proof:   Suppose that $(F|L, v)$ is henselian rational: $F^h = L(x)^h$ with $x \in F^h$. Since $(L, v)$ is henselian defectless, it is algebraically maximal. Hence in view of Theorem 3 of [KA], any pseudo–convergent sequence in $(L, v)$ without a limit in $(L, v)$ must be of transcendental type. Note that $x \notin L$ since otherwise, $F = L$. By Theorem 1 of [KA], $x$ is the limit of a pseudo–convergent sequence in $(L, v)$ without a limit in $(L, v)$, which is consequently of transcendental type. Hence by Lemma 24, $(L, v)$ is existentially closed in $(F, v)$, contradicting the fact that $L$ is not even existentially closed in $F$. This proves that $(F|L, v)$ is not henselian rational. The same argument holds with $w$ in the place of $v$. □

The function field $F$ that we have constructed shows the following symmetry between a generating Artin-Schreier extension and a generating purely inseparable extension of degree $p$. On the one hand, we have the Artin-Schreier extension
$$L(x_0, x_1)|L(x_1)$$
given by
$$x_0^p - x_0 \;=\; x - tx_1^p \;. \tag{30}$$
On the other hand we have the purely inseparable extension
$$L(x_0, x_1)|L(x_0)$$
given by
$$x_1^p \;=\; \frac{1}{t}(-x_0^p + x_0 + x) \;.$$

¿From equation (30) it is immediately clear that the function field $L(x_0, x_1)$ becomes rational after a constant field extension by $t^{1/p}$; namely

$$F(t^{1/p}) = L(t^{1/p})(x_0 + t^{1/p} x_1) .$$

This shows that the base field $L$, not being existentially closed in the function field $F$, becomes existentially closed in the function field after a finite purely inseparable constant extension, although this extension is linearly disjoint from $F|L$.

In our above example there exists also a separable constant extension $L'|L$ of degree $p$ such that $(F.L')^h$ is henselian rational. To show this, we take a constant $d \in L$ and an element $a$ in the algebraic closure of $L$ satisfying

$$t = a^p - da ,$$

and we put $L' = L(a)$. If we choose $d$ with a sufficiently high value, then we will have that $vdax_1^p > 0$. From this we deduce by Hensel's Lemma that there is an element $b \in L'(x_1)^h$ such that $b^p - b = -dax_1^p$. If we put $z = x_0 + ax_1 + b \in L'(x_0, x_1)^h$, we get that

$$z^p - z = x - tx_1^p + a^p x_1^p - ax_1 - dax_1^p = x - ax_1 + (a^p - da - t)x_1^p = x - ax_1 ,$$

which shows that

$$x_1 \in L'(z) .$$

This in turn yields that $b \in L'(z)^h$ and consequently,

$$x_0 = z - ax_1 - b \in L'(z)^h .$$

Altogether, we have proved that

$$L'(x_0, x_1)^h = L'(z)^h$$

is henselian rational.

It can be shown that extension (30) could not be immediate if $(L, v)$ resp. $(L, w)$ would satisfy property (5). This generates some hope that the crucial henselian rationality can also be proved for immediate function fields over base fields which are not tame but satisfy (PDOA).

**Modifications of our construction.**

**Modification 1:** We take $m = p - 1$ disjoint progressions $(q_{i,j})_{j \in \mathbb{N}}$, $1 \leq i \leq m$, of prime numbers $q_{i,j}$ such that for every $i$,

$$p^{j+1} < q_{i,j} \quad \text{for all } j \in \mathbb{N} . \tag{31}$$

Assume that $x_i$ are limits of the pseudo–convergent sequences

$$\left( \sum_{j=1}^{k} s^{-1/q_{i,j}} \right)_{k \in \mathbb{N}} , \quad 1 \leq i \leq m . \tag{32}$$

Then again by use of Hensel's Lemma, one shows that

$$v_s K(s, x_1, \ldots, x_k) \;=\; \sum_{i=1}^{k} \sum_{j \in \mathbb{N}} \frac{1}{q_{i,j}} \mathbb{Z}$$

for every $k \leq m$. Since the progressions are infinite and disjoint, the groups

$$v_s K(s, x_1, \ldots, x_k)/v_s K(s, x_1, \ldots, x_{k-1})$$

are not finitely generated. This proves that $x_k$ is transcendental over $K(s, x_1, \ldots, x_{k-1})$. Hence, $x_1, \ldots, x_m$ are algebraically independent over $K(s)$. Since $L_1|K(s)$ is algebraic, they are also algebraically independent over $L_1$.

We replace (17) by

$$\xi_1 = s^{-1/p} \quad \text{and} \quad \xi_{j+1} = (\xi_j - \sum_{i=1}^{m} c_i s^{-p/q_{i,j}})^{1/p} \tag{33}$$

and (24) by

$$\tilde{x}_0 := x_0 - \sum_{j=1}^{k} \xi_j \quad \text{and} \quad \tilde{x}_i := x_i - \sum_{j=1}^{k-1} s^{-1/q_{i,j}} \quad \text{for } 1 \leq i \leq m \,. \tag{34}$$

A straightforward adaptation of our arguments then shows that if $(L'|L, v)$ is an extension such that $L'v = Lv$ and $x_0, x_1, \ldots, x_m, y \in L'$ satisfy (21), then $x_1, \ldots, x_m$ are algebraically independent over $L$. Hence if $(L, v)$ is constructed with these modifications, then we obtain in addition to the assertions of Theorem 3: *If $(L'|L, v)$ is an extension such that $L'v = Lv$ and (PDOA) holds in $(L', v)$, then* trdeg $L'|L \geq p - 1$.

**Modification 2:** Instead of (5) we consider

$$\forall X \exists Y \exists X_0 \ldots \exists X_{p^n-1} \; X = Y + X_0^{p^n} - X_0 + tX_1^{p^n} + \ldots + t^{p-1} X_{p^n-1}^{p^n} \wedge \mathcal{O}(Y) \,. \tag{35}$$

We replace $m = p - 1$ by $m = p^n - 1$ and choose our prime numbers $q_{i,j}$ such that for $1 \leq i \leq p^n - 1$,

$$p^{j+n} \;<\; q_{i,j} \quad \text{for all } j \in \mathbb{N} \,. \tag{36}$$

We replace (17) by

$$\xi_1 = s^{-1/p^n} \quad \text{and} \quad \xi_{j+1} = (\xi_j - \sum_{i=1}^{m} c_i s^{-p^n/q_{i,j}})^{1/p^n} \,. \tag{37}$$

We note that by a straightforward adaptation of the proof of Lemma 4 one shows that under the assumptions of that lemma, $\mathrm{PD}(\wp(X), tX^p, \ldots, t^{p-1}X^p)$ holds. We leave the further details to the reader. If $(L, v)$ is constructed with these modifications, then we obtain in addition to the assertions of Theorem 3: *If $(L'|L, v)$ is an extension such that $L'v = Lv$ and (PDOA) holds in $(L', v)$, then* trdeg $L'|L \geq p^n - 1$.

**Modification 3:** We work in $K((s^{\mathbb{R}}))$ instead of $K((s^{\mathbb{R}}))$. We take any set of $\mathbb{Q}$-linearly independent positive real numbers $r_\ell$, $\ell \in I$ and set $s_\ell := s^{r_\ell}$. Now we do our original construction simultaneously for every $s_\ell$. That is, we extend $(K(s_\ell \mid \ell \in I), v_s)$ to a henselian defectless field $(L, w)$ with value group $wL = \sum_{\ell \in I} \mathbb{Q}r_\ell$. With this modification, we obtain: *If $(L'|L, v)$ is an extension such that $L'v = Lv$ and (PDOA) holds in $(L, v)$, then* $\operatorname{trdeg} L'|L \geq (p-1) \cdot |I|$. Hence if $I$ is infinite, then the transcendence degree will be infinite. However, the modified extension $L|K$ will also be of transcendence degree $|I|$; in particular, $L$ will not be a field of finite transcendence degree over its prime field. So we ask:

**Can $(L, v)$ be constructed in such a way that all assertions of Theorem 3 hold and that $\operatorname{trdeg} L'|L$ is infinite for every extension $(L'|L, v)$ for which $L'v = Lv$ and (PDOA) holds in $(L', v)$?**

## 5 Appendix: Images of polynomials in valued fields

By a **generalized ball** in a valued field $(K, v)$ we mean a union over any nest of balls. Note that $B$ is a generalized ball in $(K, v)$ if and only if for all $a, b \in B$ we have that $B_{v(a-b)}(a) \subset B$.

Generalized balls have the same property as ordinary balls: if the intersection of two of them is nonempty, then they are comparable by inclusion. Hence we can define a **nest of generalized balls** in the same way as a nest of balls. For every nest $\mathbf{B}$ of generalized balls there is a nest $\overline{\mathbf{B}}$ of closed balls such that for every ball $B \in \mathbf{B}$ there is some ball $\overline{B} \in \overline{\mathbf{B}}$ with $\overline{B} \subset B$. This can be found as follows. If $B \in \mathbf{B}$ is not the smallest ball in $\mathbf{B}$, then there is some $a \in B$ which is not contained in any of the smaller balls. We pick some smaller ball $B_0 \in \mathbf{B}$ and an element $a_0 \in B_0$. Then the closed ball $\overline{B} := B_{v(a-a_0)}(a)$ is contained in $B$ and contains all smaller balls of $\mathbf{B}$. If $B \in \mathbf{B}$ is the smallest ball in $\mathbf{B}$, then we pick some $a \in B$. If $B$ is a singleton, then $B = B_\infty(a)$ is already a closed ball. Otherwise, we pick some $b \in B$ different from $a$ and set $\overline{B} := B_{v(a-b)}(a) \subset B$. By construction, the so-obtained set $\{\overline{B} \mid B \in \mathbf{B}\}$ is a nest of closed balls, and we have that $\bigcap \overline{\mathbf{B}} \subseteq \bigcap \mathbf{B}$. Equality holds if B contains no smallest ball. Hence if $(K, v)$ is spherically complete, then every nest of generalized balls in $(K, v)$ has a nonempty intersection.

Take any valued field $(K, v)$ and let $f$ be any map of $K$ into itself. Assume that

$$\left.\begin{array}{l}\text{for every nest } \mathbf{B}' \text{ of closed balls in the image } f(K) \text{ of } f \\ \text{there is a nest } \mathbf{B} \text{ of generalized balls in } K \text{ such that:} \\ \text{for every } B' \in \mathbf{B}' \text{ there is } B \in \mathbf{B} \text{ satisfying that } f(B) \subseteq B'.\end{array}\right\} \quad (38)$$

Then $f(a) \in \bigcap_{B' \in \mathbf{B}'} B'$ for every element $a \in \bigcap_{B \in \mathbf{B}} B$. Hence,

**Lemma 29** *If $(K, v)$ is spherically complete and $f$ has property (38), then $(f(K), v)$ is spherically complete.*

To see that all polynomials $f$, viewed as maps on valued fields, have the property (38), we just have to describe the preimages of closed balls under $f$. Beforehand, we note the following. If $(L, v)$ is any extension of $(K, v)$ and $B$ is a generalized ball in $(L, v)$, then $B \cap K$ is a generalized ball in $(K, v)$, or empty. Hence if $\mathbf{C}$ is any collection of generalized balls in $(L, v)$, then

$$\mathbf{C} \cap K := \{B \cap K \mid B \in \mathbf{C} \text{ such that } B \cap K \neq \emptyset\}$$

is a collection of generalized balls in $(K,v)$. If $\mathbf{B}$ is a nest of generalized balls in $(L,v)$, then $\mathbf{B}\cap K$ is a nest of generalized balls in $(K,v)$, if nonempty.

**Lemma 30** *Let $f\in K[X]$ be of degree $n$ and $B'$ a generalized ball in $(K,v)$. Then the preimage $f^{-1}(B')$, if nonempty, is the disjoint union of at most $n$ uniquely determined generalized balls in $(K,v)$. If in addition $K$ is algebraically closed and $c$ is any element of $B'$, then each of these balls contains at least one of the roots of the polynomial $f-c$.*

Proof: First, we show our assertion for the case of $K$ algebraically closed. Let $B'$ be a generalized ball in $(K,v)$. Let $d$ be the leading coefficient and $a_1,\ldots,a_n$ the roots of the polynomial $f-c$. Then

$$v(f(a)-c) \;=\; v(f-c)(a) \;=\; v\left(d\prod_{i=1}^{n}(a-a_i)\right) \;=\; vd + \sum_{i=1}^{n} v(a-a_i)\;.$$

Assume that $a\in f^{-1}(B')$ and pick $k\in\{1,\ldots,n\}$ such that $v(a-a_k)=\max_i v(a-a_i)$. Now if $v(b-a_k)\geq v(a-a_k)$ then

$$v(b-a_i) \;\geq\; \min\{v(b-a_k),v(a_k-a),v(a-a_i)\} \;=\; v(a-a_i)$$

for all $i$. It follows that

$$v(f(b)-c) \;=\; vd + \sum_{i=1}^{n} v(b-a_i) \;\geq\; vd + \sum_{i=1}^{n} v(a-a_i) \;=\; v(f(a)-c)\;,$$

which yields that $f(b)\in B'$. That is, $b\in f^{-1}(B')$. Therefore,

$$B_{v(a-a_k)}(a_k) \subset f^{-1}(B')\;.$$

We set

$$B_k \;:=\; \bigcup\{B_{v(a-a_k)}(a_k) \mid a\in f^{-1}(B') \text{ and } v(a-a_k)=\max_i v(a-a_i)\}\;, \tag{39}$$

which is a generalized ball in $(K,v)$, or empty. We conclude that

$$f^{-1}(B') \;=\; \bigcup_{1\leq k\leq n} B_k\;.$$

If $B_j\cap B_k\neq\emptyset$, then one of them contains the other, say, $B_j\subseteq B_k$. Thus, to obtain a disjoint union, we can just omit $B_j$. This proves the first assertion.

The fact that generalized balls with nonempty intersection are comparable by inclusion can easily be used to show that the generalized balls in a disjoint union are uniquely determined. Since $c\in B'$ was arbitrarily chosen, this also yields our second assertion.

Now assume that $K$ is not algebraically closed. We choose any extension of $v$ to $\tilde{K}$. We associate a generalized ball $\tilde{B}'$ to $B'$ in the following way: we pick some $a\in B'$ and set

$$\tilde{B}' \;:=\; \bigcup\{B_\alpha(a,\tilde{K}) \mid \alpha\in vK\cup\{\infty\} \text{ such that } B_\alpha(a,K)\subset B'\}\;. \tag{40}$$

Then $\tilde{B}' \cap K = B'$. By what we have proved, $f^{-1}(\tilde{B}')$ is a disjoint union of at most $n$ generalized balls in $(\tilde{K}, v)$. Hence, $f^{-1}(B') = f^{-1}(\tilde{B}') \cap K$ is the disjoint union of at most $n$ generalized balls in $(K, v)$, or empty. □

For the sake of completeness, let us mention that $B_j \cap B_k \ne \emptyset$ implies that $B_j = B_k$. To see this, assume that $B_j \subseteq B_k$. Then in particular, $a_j \in B_k$. That is, there is some $a \in B_k$ such that $v(a - a_k) = \max_i v(a - a_i)$ and $a_j \in B_{v(a-a_k)}(a_k)$. But then, $v(a - a_j) \le v(a - a_k)$ by the former and $v(a - a_j) \ge v(a - a_k)$ by the latter. Hence, $v(a - a_j) = v(a - a_k)$ and the ball $B_{v(a-a_k)}(a_k) = B_{v(a-a_j)}(a_j)$ is contained in both $B_j$ and $B_k$. If $a' \in B_k$ is such that $v(a' - a_k) = \max_i v(a' - a_i)$ and $B_{v(a-a_k)}(a_k) \subseteq B_{v(a'-a_k)}(a_k)$, then $a_j \in B_{v(a'-a_k)}(a_k)$ and the same argument shows that $B_{v(a'-a_k)}(a_k) \subseteq B_j$. It follows that $B_k \subseteq B_j$, so $B_j = B_k$.

In [K2], we show that if $B'$ is a closed ball in $(K, v)$, then $f^{-1}(B')$, if nonempty, is the disjoint union of at most $n$ uniquely determined closed balls in $(K, v)$, provided that $K$ is algebraically closed (or, more generally, that $vK$ is divisible).

Using the well-known theorem that the inverse limit of an inverse system of nonempty compact hausdorff spaces is nonempty, one can prove the following:

**Lemma 31** *Let $f$ be a map on an ultrametric space $Y$. Suppose that there is a natural number $n$ such that the preimage under $f$ of every generalized ball is the union of at most $n$ generalized balls. Then $f$ has property (38).*

It follows that every polynomial $f$ on a valued field has property (38). But here, we wish to deduce this fact in a different way by using a stronger property of polynomial maps. For any nest $\mathbf{B}'$ in $(K, v)$ we define $f^{-1}(\mathbf{B}')$ to be the collection of all generalized balls which for some $B' \in \mathbf{B}'$ appear in the disjoint union representation of $f^{-1}(B')$ given by Lemma 30.

If $(L, v)$ is any extension of $(K, v)$ and $\mathbf{C}$ is a union of at most $n$ nests of generalized balls in $(L, v)$, then $\mathbf{C} \cap K$ is a union of at most $n$ nests of generalized balls in $(K, v)$. We leave it as a simple exercise to the reader to prove that a collection $\mathbf{C}$ of generalized balls is the union of at most $n$ nests of generalized balls if and only if the maximal number of pairwise disjoint generalized balls in $\mathbf{C}$ is at most $n$.

**Lemma 32** *For any nest $\mathbf{B}'$ of generalized balls in $(K, v)$ and every polynomial $f \in K[X]$ of degree $n$, $f^{-1}(\mathbf{B}')$ is a union of at most $n$ nests of generalized balls.*

Proof:  First, we show our assertion for the case of $K$ algebraically closed. For every choice of pairwise disjoint generalized balls $B^{(1)}, \ldots, B^{(k)} \in f^{-1}(\mathbf{B}')$ we have to show that $k \le n$. There are finitely many balls $B'_1, \ldots, B'_m \in \mathbf{B}'$ such that every $B^{(i)}$ appears in the disjoint union representation of $f^{-1}(B'_j)$ for some $j$. Since $B'_1, \ldots, B'_m$ are linearly ordered by inclusion, we can choose some $c$ in their intersection. Let $a_1, \ldots, a_n$ be the roots of $f - c$. Then by Lemma 30, every $B^{(i)}$ contains at least one $a_j$. Hence, $k \le n$.

Now assume that $K$ is not algebraically closed. We choose any extension of $v$ to $\tilde{K}$. To our nest $\mathbf{B}'$ we associate a nest $\tilde{\mathbf{B}}'$ of generalized balls in $(\tilde{K}, v)$ in the following way:

$$\tilde{\mathbf{B}}' := \{\tilde{B}' \mid B' \in \mathbf{B}'\},$$

where $\tilde{B}'$ is defined as in (40). Then $\tilde{\mathbf{B}}' \cap K = \mathbf{B}'$. By what we have proved above, $f^{-1}(\tilde{\mathbf{B}}')$ is a union of at most $n$ nests of generalized balls in $(\tilde{K}, v)$. Hence, $f^{-1}(\mathbf{B}') = f^{-1}(\tilde{\mathbf{B}}') \cap K$ is a union of at most $n$ nests of generalized balls in $(K, v)$. □

Among the finitely many nest of generalized balls that $f^{-1}(\mathbf{B}')$ consist of, at least one nest $\mathbf{B}$ must satisfy (38). This proves:

**Lemma 33** *Take any valued field $(K, v)$. Then all polynomials $f \in K[X]$ have the property (38).*

By virtue of Lemma 29, this proves Lemma 12.

**Proof of Lemma 13:** Suppose that $f \in K[X]$ is a polynomial such that $f(K)$ is not an optimal approximation subset of $(K, v)$. Take some $c \in K$ such that $\{v(c-x) \mid x \in f(K)\}$ does not have a maximum. Then the nest

$$\mathbf{B}' := \{B_{v(c-x)}(c, K) \cap f(K) \mid x \in f(K)\}$$

of balls in $f(K)$ does not contain a smallest ball, and $\bigcap \mathbf{B}'$ is empty. By Lemma 33 we can choose a nest $\mathbf{B}$ of generalized balls in $(K, v)$ such that (38) holds. Then also $\bigcap \mathbf{B}$ is empty. Therefore, $\mathbf{B}$ does not contain a smallest ball, and thus its coinitiality is a limit ordinal, say, $\lambda$. We pick a coinitial chain $(B_\nu)_{\nu < \lambda}$ of generalized balls in $\mathbf{B}$. Then for every $\nu < \lambda$ we pick some $b_\nu \in B_\nu \setminus B_{\nu+1}$. We leave it as an exercise to the reader to show that $(b_\nu)_{\nu < \lambda}$ is a pseudo–convergent sequence. If it would have a limit $b$ in $K$, then $b$ would lie in $\bigcap_{\nu < \lambda} B_\nu = \bigcap \mathbf{B}$. But the latter intersection is empty, so the sequence has no limit in $K$.

Suppose that there were some $\nu_0 < \lambda$ such that for all $\nu \geq \nu_0$, the value $v(f(b_\nu) - c)$ is constant, say, equal to $\beta$. As $\mathbf{B}'$ contains the nonempty ball $B_\alpha(c, K) \cap f(K)$ but does not contain a smallest ball, there must be some ball in $\mathbf{B}'$ of radius $> \beta$. This ball contains none of the elements $f(b_\nu)$, $\nu \geq \nu_0$. But this contradicts the fact that B satisfies (38) and $(B_\nu)_{\nu < \lambda}$ is a coinitial sequence in $\mathbf{B}$, with $b_\nu \in B_\nu$. We thus find that the sequence $(b_\nu)_{\nu < \lambda}$ does not fix the value of the polynomial $f - c$. This shows that $(b_\nu)_{\nu < \lambda}$ is of algebraic type. Take $g \in K[X]$ to be a monic polynomial of minimal degree such that $(b_\nu)_{\nu < \lambda}$ does not fix the value of $g$. If $g$ were linear, say equal to $X - b$, then $b$ would be a limit of $(b_\nu)_{\nu < \lambda}$. But $(b_\nu)_{\nu < \lambda}$ has no limit in $K$, so $\deg g > 1$. Now if $b$ is any root of $g$, then by Theorem 3 of [KA], there is an immediate extension of $v$ from $K$ to the proper algebraic extension $K(b)$, showing that $(K, v)$ is not algebraically maximal. □

# References

[AK]   Ax, J. – Kochen, S.: *Diophantine problems over local fields I, II*, Amer. Journ. Math. **87** (1965), 605–630, 631–648

[D]    Delon, F.: *Quelques propriétés des corps valués en théories des modèles*, Thèse Paris VII (1981)

[E]    Ershov, Yu. L.: *On the elementary theory of maximal valued fields I, II, III* (in Russian), Algebra i Logika **4**:3 (1965), 31–70, **5**:1 (1966), 5–40, **6**:3 (1967), 31–38

[G]        Gravett, K. A. H.: *Note on a result of Krull*, Cambridge Philos. Soc. Proc. **52** (1956), 379

[KA]       Kaplansky, I.: *Maximal fields with valuations I*, Duke Math. Journ. **9** (1942), 303–321

[KR]       Krull, W.: *Allgemeine Bewertungstheorie*, J. reine angew. Math. **167** (1931), 160–196

[K1]       Kuhlmann, F.-V.: *Henselian function fields and tame fields*, preprint (extended version of Ph.D. thesis), Heidelberg (1990)

[K2]       Kuhlmann, F.-V.: *Valuation theory of fields, abelian groups and modules*, preprint, Heidelberg (1996), to appear in the "Algebra, Logic and Applications" series (Gordon and Breach), eds. A. Macintyre and R. Göbel

[K3]       Kuhlmann, F.-V.: *On local uniformization in arbitrary characteristic*, The Fields Institute Preprint Series, Toronto (July 1997)

[K4]       Kuhlmann, F.-V.: *A theorem about maps on spherically complete ultrametric spaces*, preprint, Urbana (November 1997)

[K5]       Kuhlmann, F.-V.: *Valuation theoretic and model theoretic aspects of local uniformization*, submitted for the Proceedings of the Blowup Tirol Conference 1997, Saskatoon (1998)

[K6]       Kuhlmann, F.-V.: *On places of algebraic function fields in arbitrary characteristic*, preprint, Toronto/Saskatoon (1997)

[K7]       Kuhlmann, F.-V.: *The model theory of tame valued fields*, in preparation, Saskatoon (1998)

[K8]       Kuhlmann, F.-V.: *On local uniformization in arbitrary characteristic I*, in preparation, Saskatoon (1998) (improved version of [26])

[K9]       Kuhlmann, F.-V.: *On local uniformization in arbitrary characteristic II*, in preparation, Saskatoon (1998)

[KP]       Kuhlmann, F.-V. – Prestel, A.: *On places of algebraic function fields*, J. reine angew. Math. **353** (1984), 182–195

[KPR]      Kuhlmann, F.-V. – Pank, M. – Roquette, P.: *Immediate and purely wild extensions of valued fields*, manuscripta math. **55** (1986), 39–67

[L]        Lang, S.: *Algebra*, Addison-Wesley, New York (1965)

[O]        Ore, O.: *On a special class of polynomials*, Trans. Amer. Math. Soc. **35** (1933), 559–584

[R]        Ribenboim, P.: *Théorie des valuations*, Les Presses de l'Université de Montréal, Montréal, 1st ed. (1964), 2nd ed. (1968)

[W1]       Whaples, G.: *Additive polynomials*, Duke Math. Journ. **21** (1954), 55–65

[W2]       Whaples, G.: *Galois cohomology of additive polynomials and n-th power mappings of fields*, Duke Math. Journ. **24** (1957), 143–150

[Z]        Ziegler, M.: Die elementare Theorie der henselschen Körper, *Inaugural Dissertation*, Köln (1972)

Department of Mathematics and Statistics, University of Saskatchewan,
106 Wiggins Road, Saskatoon, Saskatchewan, Canada S7N 5E6
email: fvk@math.usask.ca   —   home page: http://math.usask.ca/~fvk/index.html